\documentclass[letterpaper]{article}

\usepackage{packages-styles/general-math}
\usepackage{packages-styles/field-probabilitytheory}
\usepackage{packages-styles/field-algorithmic}
\usepackage{packages-styles/general-scientific}
\newcommand\Language{english}
\usepackage{packages-styles/general-text}
\crefname{assumption}{Assumption}{Assumptions}

\usepgfplotslibrary{external} 
\usepackage{xcolor}
\usepackage{arxiv}

\hypersetup{
pdftitle={A Principle for Global Optimization with Gradients},
pdfsubject={Mathematical Optimization},
pdfauthor={Nils M{\"u}ller},
pdfkeywords={Global Optimization, Robust Optimization, Continuous Optimization, Mathematical Optimization, Simulation-Based Optimization}
}

\newcommand{\NN}{\mathcal{N}}
\Crefname{property}{Property}{Properties}
\setlist[description]{font=\rmfamily}

\crefname{algorithm}{Algorithm}{Algorithms}
\Crefname{algorithm}{Alg.}{Algs.}

\title{A Principle for Global Optimization with Gradients}

\author{
Nils M{\"u}ller\\[0.3em]
Max Planck Institute for Software Systems\\
Saarland Informatics Campus\\[0.3em]
Max Planck Institute for Intelligent Systems\\[0.3em]
Germany\\[0.3em]
\href{mailto:nmueller@mpi-sws.org}{\texttt{nmueller@mpi-sws.org}}
}

\begin{document}
\maketitle

\begin{abstract}
    This work demonstrates the utility of gradients for the global optimization of certain differentiable functions with many suboptimal local minima. To this end, a principle for generating search directions from non-local quadratic approximants based on gradients of the objective function is analyzed. Experiments measure the quality of non-local search directions as well as the performance of a proposed simplistic algorithm, of the  covariance matrix adaptation evolution strategy (CMA-ES), and of a randomly reinitialized Broyden-Fletcher-Goldfarb-Shanno (BFGS) method.
\end{abstract}


\section{Introduction}
This work motivates the use of gradients for solving the optimization problem $\min_{x \in \R^n} f(x)$, where the differentiable function $f: \R^n \to \R$ is assumed to have many suboptimal local minima but possesses (unknown) global structure.\\
Optimization methods based on the (quasi-)Newton method, see, e.g., \cite[Chapter~8]{nocedal}, undesirably converge when applied to the minimization of functions with many suboptimal local minima. Contrary, in a formal setting, non-local information (as opposed to global) provides search directions for which iterative optimizers converge to the global minimum of certain functions with many suboptimal local minima \cite[Theorem~4.1~\&~Theorem~5.2]{mueller}. Considering in addition that non-local information based on objective evaluations has been successfully used in many practical optimization methods for decades, this work develops a (quasi-)Newton method that approximates search directions from non-local gradient information of practically realistic evaluation counts. To achieve this, the search direction, commonly based on a quadratic Taylor-approximant in (quasi-)Newton methods, is replaced with a non-local generalization, i.e., \texttt{line 3} of \cref{alg:nlqn}. \cref{alg:nlqn} is a simplistic method with the purpose of demonstrating the usefulness of gradients beyond local optimization.\\

The interest in investigating the utility of gradients for the global optimization of functions with many suboptimal local minima lies in the hope of generalizing their established success from local optimization. Computer simulation, in particular in the context of the adjoint method \cite{lions,pironneau,cea} and automatic differentiation \cite{wengert,linnainmaa}, can provide gradient information on challenging simulation-based objectives from science and engineering \cite{assad,tromp,othmer}. Therefore, associated optimization methods based on non-local evaluations of the objective gradients may have a large impact on a wide range of real-world challenges.

\begin{algorithm}[h]
\caption{Non-local quasi-Newton method\protect\footnotemark[1] for optimization of differentiable $f: \R^n \to \R$%
}
\begin{algorithmic}[1]
    \INPUT Continuously differentiable function $f: \R^n \to \R$, and its gradient $\nabla f: \R^n \to \R^n$;
    \SUBINPUT initial point $x_0 \in \R^n$;
    \SUBINPUT initial scaling $\sigma_0 \in \R_{>0}$;
    \SUBINPUT sample size $k \in \N$;
    \SUBINPUT a probability measure $\Pr_k$ on $\R^{n \times k}$;
    \SUBINPUT non-local line search method $\mathrm{linesearch}(\cdot)$;
    \SUBINPUT scaling adaption method $\mathrm{scaling}(\cdot)$;
    \SUBINPUT maximum iteration count $C \in \N$.
    \OUTPUT Element $x_C \in \R^n$ with \enquote{low} function value $f(x_C)$.
    \INIT $q_{A, b}(x) := \langle x, (A+A^T) x \rangle + b^T x$ for all $x \in \R^n$, $A \in \R^{n \times n}$, $b \in \R^n$; $t := 0$.
    \AlgBlankLine
    \STATE \texttt{sample} $z_1, \dots, z_k \in \R^n$ according to $\Pr_k$
    \algcomment{independent of any other time-step}
    \STATE \texttt{compute} $\nabla f(x_t +\sigma_t z_1), \dots, \nabla f(x_t + \sigma_t z_k) \in \R^n$
    \algcomment{sample a neighborhood of $x_t$}
    \STATE \texttt{solve for} $\Delta x_t \in \R^n$        \algcomment{determine non-local Newton direction\protect\footnotemark[2]}
    \[\begin{cases} \Delta x_t \in \displaystyle\argmin_{x \in \R^n} \,q_{A_t, b_t}(x)\\
    (A_t, b_t) \in \displaystyle\argmin_{A \in \R^{n\times n},\, b \in \R^n }\sum_{j=1}^k \norm{\nabla q_{A,b}(z_j) - \nabla f(x_t + \sigma_t z_j)}^2
    \end{cases}\]
    \STATE $x_{t+1} := \mathrm{linesearch}(f, \Delta x_t, -b_t, x_t)$
    \algcomment{non-local linesearch based on $\Delta x_t$ and $-b_t$}
    \STATE $\sigma_{t+1} := \mathrm{scaling}(\sigma_0, \sigma_t, x_{t+1}-x_t)$
    \algcomment{adapt scaling}
    \STATE $t \leftarrow t + 1$
    \STATE \texttt{if} $t < C$ \texttt{go to} $\mathrm{line}$ 1
    \algcomment{until the budget is used}
    \STATE \textbf{return} $x_{C}$
  \end{algorithmic}
\label{alg:nlqn}
\end{algorithm}
\refstepcounter{footnote}%
\footnotetext{An implementation in the Python programming language can be found at \href{https://github.com/NiMlr/pynlqn}{\texttt{https://github.com/NiMlr/pynlqn}}.}%
\refstepcounter{footnote}%
\footnotetext{In case a minimum of $q_{A_t, b_t}$ is not attained on $\R^n$, consider a trust region, as specified in \cref{rem:trust}.}%

\noindent\textbf{Related Work.} Non-local operators for optimization have been shown to impose favorable structure on objective approximations \cite{more,loog,mobahi,mueller} and can be asymptotically consistent with local operators \cite{nagaraj}, i.e., when the non-local kernel converges to the Dirac-measure. \cite{omeradzic} point out the open question of sample-efficient approximations of non-local operators, while \cite{berahas} indicate the superiority of first-order interpolations from function evaluations over empirical gradient averages. \cite{zhang} propose smoothing with a non-local kernel with low-dimensional support. \cite{maggiar,berahas2} propose a similar approach to the one of this work using quadratic approximation from function instead of gradient evaluations.\\
A (quasi-)Newton method roughly similar to \cref{alg:nlqn} for the different goal of locally solving systems of equations based on a least-squares principle has been described by \cite{haelterman}.\\
Practical methods employing non-local operators of function evaluations have been widely used, a few examples include \cite{beyer,salimans,hansen,ollivier,schwefel,rechenberg,storn}.\\
Specifically, attempts have been made to design methods that combine the advantages of non-local operators of function evaluations and gradient methods, however, in contrast to this work they do not make use of explicit gradient information \cite{yalcin,salomon,arnold}.\\
Under the assumption of a~priori~knowledge of non-local information about the objective function at hand, even local gradient methods have been developed for global optimization \cite{griewank}.\\
The proposed principle can also be considered a non-local gradient-based trust region method. Trust region methods are commonly restricted to finding local-optima \cite[Chapter~3.2]{conn}, while non-local extensions have been restricted to using function evaluations, i.e., a black-box setting \cite{addis,diouane,cheng}.\\

\noindent\textbf{Outline.} In \cref{thm:nsc,cor:nsc}, necessary and sufficient conditions that let us compute the proposed non-local search direction based on a solution to a Lyapunov-type equation are developed. Under the assumption of infinite sampling, it is then proven in \cref{thm:lanc} that the solution of \texttt{line 3} of \cref{alg:nlqn} is a consistent approximator of the optimal search direction on a quadratic that is disturbed by a function with bounded first derivative. The theoretical analysis is concluded in \cref{thm:resb} by a non-local residual bound for Rastrigin-type objective functions which present an important model for functions with many suboptimal local minima.\\
In a first experiment in \cref{subsec:exp1}, gradient-based search directions for a Rastrigin-type objective function are compared. Further, \cref{alg:nlqn}, the \emph{Covariance matrix adaptation evolution strategy (CMA-ES)}, and a randomly reinitialized \emph{Broyden-Fletcher-Goldfarb-Shanno method} are benchmarked on selected functions with many local minima in \cref{subsec:exp2}. Lastly, it is shown in \cref{subsec:exp3} that \cref{alg:nlqn} solves \emph{Problem 4 of the SIAM News: A Hundred-dollar, Hundred-digit Challenge}. The work concludes with a discussion of the results and promising future work in \cref{sec:discussion}.

\section{Analysis}
\label{sec:analysis}

Initially, necessary and sufficient conditions that let us efficiently compute the approximant $q_{A_t, b_t}$ in \cref{alg:nlqn} are developed. To this end, consider the following \cref{thm:nsc}.
\begin{theorem}[Necessary and sufficient conditions for ($A_t, b_t$)]
\label{thm:nsc}
In the setting of \cref{alg:nlqn}, $A_t \in \R^{n\times n}$ and $b_t \in \R^n$ are the parameters of a best quadratic approximation in the sense of \texttt{line 3} of \cref{alg:nlqn} if and only if
\[
\begin{cases}
    (A_t + A_t^T) (Z- \overline{Z})Z^T + Z(Z -  \overline{Z})^T (A_t + A_t^T) = (G-\overline{G})Z^T + Z(G-\overline{G})^T 
    \\
    b_t = \overline{g}- (A_t+A_t^T) \overline{z}
    \,,
\end{cases}
\]
where
\begin{itemize}
    \item $G \in \R^{n \times k}$ by $G_{\cdot, j} := \nabla f(x_t + \sigma_t z_j)$ for all $j \in \N_{\leq k}$,
    \item $Z \in \R^{n \times k}$ by $Z_{\cdot, j} := 2\sigma_t z_j$ for all $j \in \N_{\leq k}$,
    \item $\overline{z} := (1/k)\sum_{j=1}^k Z_{\cdot, j}$,
    \item $\overline{g} := (1/k)\sum_{j=1}^k G_{\cdot, j}$,
    \item $\overline{Z} \in \R^{n \times k}$ by $\overline{Z}_{\cdot, j} := \overline{z}$ for all $j \in \N_{\leq k}$, and
    \item $\overline{G} \in \R^{n \times k}$ by $\overline{G}_{\cdot, j} := \overline{g}$ for all $j \in \N_{\leq k}$.
\end{itemize}
\end{theorem}
\begin{proof}
    In the setting of \cref{alg:nlqn},
    \begin{itemize}
        \item one has $\nabla q_{A,b}(\sigma_t z_j) = (A + A^T)2\sigma_t z_j + b$
        \item define $B \in \R^{n \times k}$ by $B_{\cdot, j} := b$ for all $j \in \N_{\leq k}$, and
        \item let $\norm{\cdot}_{\mathrm{F}}$ be the Frobenius norm on $\R^{n \times k}$.
    \end{itemize}
    \emph{i.~Representation of the objective.} One observes that
    \begin{align*}
        \sum_{j=1}^k \norm{\nabla q_{A,b}(z_j) - \nabla f(x_t + \sigma_t z_j)}^2
        &=
        \sum_{j=1}^k \norm{(A+A^T) Z_{\cdot, j} + B_{\cdot, j} - G_{\cdot, j}}^2
    \tag*{\small{(definitions $G, Z, B$)}}
        \\
        &=
        \norm{(A+A^T) Z + B - G}_{\mathrm{F}}^2
    \tag*{\small{(definition of $\norm{\cdot}_{\mathrm{F}}$)}}
        \\
        &=: W(A, b) \,.
    \end{align*}
    \emph{ii.~First-order conditions.} For all $A \in \R^{n \times n}$, one has a simple first order criterion an optimum $b \in \R^n$ of $W(A, \cdot)$, which reads
    \begin{align*}
        \frac{(\partial W)(A,b)}{\partial b}
        &=
        k\big( (A+A^T) \overline{z} + b - \overline{g} \big) \overset{!}{=} 0
        \\
        &\quad\iff
        b = \overline{g}- (A+A^T) \overline{z}
        \\
        &\quad\iff
        B = \overline{G}- (A+A^T) \overline{Z} \,.
    \tag*{\small{(definitions $\overline{G}, \overline{Z}$)}}
    \end{align*}
    Further, it is known that for all $b \in \R^n$ a minimizer $A \in \R^{n \times n}$ of $W(\cdot, b)$ fulfills
    \begin{align*}
        \frac{(\partial W)(A,b)}{\partial A}
        &=
        2 \big( (A+A^T) Z + B - G \big) Z^T + 2 Z \big( (A+A^T)Z + B - G \big)^T \overset{!}{=} 0
        \\
        &\quad\iff
        (A+ A^T) Z Z^T + ZZ^T(A + A^T) = (G-B)Z^T + Z(G-B)^T
        \\
        &\quad\iff
        (A+ A^T) Z Z^T + ZZ^T(A + A^T)
        \\
        &\qquad\qquad\qquad =
        \big(G-\overline{G}+ (A+A^T) \overline{Z}\big)Z^T + Z\big(G-\overline{G}+ (A+A^T) \overline{Z}\big)^T
    \tag*{\small{(assuming $B = \overline{G}- (A+A^T) \overline{Z}$)}}
        \\
        &\quad\iff
        (A+ A^T) Z Z^T + ZZ^T(A + A^T)
        \\
        &\qquad\qquad\qquad =
        (G-\overline{G})Z^T+ (A + A^T) \overline{Z}Z^T + Z(G-\overline{G})^T+ Z \overline{Z}^T (A + A^T)
        \\
        &\quad\iff
        (A + A^T) (Z- \overline{Z})Z^T + Z(Z -  \overline{Z})^T (A + A^T) = (G-\overline{G})Z^T + Z(G-\overline{G})^T \,.
    \end{align*}
    \emph{iii.~Sufficiency argument.} As the symmetric matrices are a linear subspace of $\R^{n \times n}$, and in particular, without boundary and closed, any minimizer of $W$ must be a critical point. Further, as the problem is also convex, it can be concluded that any critical point of $W$ is also a minimizer.
\end{proof}

In fact, by simple first-order conditions, the search direction is also obtained. The second equation of the following \cref{cor:nsc} is an equation of Lyapunov-type for which a wide range of numerical methods exist.
\begin{corollary}[Necessary and sufficient conditions for $\Delta x_t$]
\label{cor:nsc}
In the setting of \cref{alg:nlqn} and \cref{thm:nsc}, $\Delta x_t$ is a well-defined minimizer of a well-defined $q_{A_t, b_t}$ if there exists a positive definite matrix $\Tilde{A}_t \in \R^{n \times n}$ such that
\begin{align*}
\begin{cases}
    (\Tilde{A}_t + \Tilde{A}_t^T) \Delta x_t  = -b_t
    \\
    \Tilde{A}_t (Z- \overline{Z})Z^T + Z(Z -  \overline{Z})^T \Tilde{A}_t = (G-\overline{G})Z^T + Z(G-\overline{G})^T \,,
\end{cases}
\end{align*}
where $b_t = \overline{g}- (A_t+A_t^T) \overline{z}$ and $A_t = (\Tilde{A}_t + \Tilde{A}_t^T)/4$.
\end{corollary}
\begin{remark}
\label{rem:trust}
    In case there is no positive definite matrix $\Tilde{A}_t \in \R^{n\times n}$ that fulfills the conditions of \cref{cor:nsc}, consider finding a minimum in a trust region $D^n := \{ x \in \R^n \mid \norm{x} \leq 1\}$, i.e.,
    \[
        \Delta x_t \in \argmin_{x \in D^n} q_{A_t, b_t}(x) \,,
    \]
    where $A_t, b_t$ are picked as specified in \cref{cor:nsc} for a possibly negative definite or indefinite matrix $\Tilde{A}_t$ that fulfills the remaining conditions.
\end{remark}
\begin{proof} Define $P := (Z- \overline{Z}) Z^T$ and $V:= (G- \overline{G})Z^T$.\\

\noindent\emph{i.~It is claimed that a matrix $A_t \in \R^{n\times n}$ that fulfills the conditions of \cref{thm:nsc} exists if and only if there exists $\Tilde{A}_t \in \R^{n \times n}$ with
    $\Tilde{A}_t P + P^T \Tilde{A}_t = V + V^T$.}
    Assuming such $\Tilde{A}_t$ exists and $A_t := (\Tilde{A}_t + \Tilde{A}^T_t)/4$, one has
    \begin{align*}
        (A_t + A_t^T) P + P^T (A_t + A_t^T)
        &=
        \big((\Tilde{A}_t + \Tilde{A}_t^T)/2\big) P + P^T (\Tilde{A}_t + \Tilde{A}_t^T)/2
    \tag*{\small{(definition and symmetry of $A_t$)}}
        \\
        &=
        (\Tilde{A}_t P + P^T \Tilde{A}_t )/2 + ( \Tilde{A}_t P + P^T \Tilde{A}_t)^T/2
    \tag*{\small{(distributivity, commutativity of addition, transpose of products)}}
        \\
        &=
        (V + V^T )/2 + ( V + V^T)^T/2
    \tag*{\small{(condition on $\Tilde{A}_t$)}}
        \\
        &=
        V + V^T \,.
    \tag*{\small{(distributivity, commutativity of addition)}}
    \end{align*}
    The converse statement is true as for $A_t \in \R^{n¸\times n}$ that fulfills the condition of \cref{thm:nsc}, the substitution $\Tilde{A}_t := A_t + A_t^T$ generates the condition that is to prove. Without loss of generality $A_t$ is symmetric---otherwise take $(A_t + A_t^T)/2$ instead of $A_t$.\\

    \noindent\emph{ii.~A first order condition on $\Delta x_t \in \argmin_{x \in \R^n} q_{A_t, b_t}(x)$ is necessary and sufficient due to convexity of $q_{A_t, b_t}$, and reads
    \[
        (\nabla q_{A_t, b_t}) (\Delta x_t) = 2(A_t + A_t^T) \Delta x_t  + b_t \overset{!}{=} 0 \,.
    \qedhere
    \]}
\end{proof}

For simplicity, the following analysis is restricted to relating the non-local approximation objective of determining a quadratic model in \texttt{line 3} of \cref{alg:nlqn} with a $k$-asymptotic setting by the following \cref{ass:gc}. The assumption models a setting, where an infinite amount of non-local gradient samples is available. Note, that $\Pr_k$ is a measure on $\R^{n \times k}$, e.g., the $k$-product of an $n$-variate normal distribution or the $k$-product of (distinct) Dirac measures in $\R^n$.

\begin{assumption}[Glivenko-Cantelli]
\label{ass:gc}
    In \cref{alg:nlqn}, let w.l.o.g.~$x_t = 0$, otherwise transform $f$. Assume that for all $\sigma \in (0, \infty)$, some probability measure $\Pr_\sigma := \Pr(\sigma^{-1} \cdot)$ on $\R^n$, probability measure $\Pr_\infty$ on $\R^\N$ with $(z_1, z_2, \dots) \sim \Pr_\infty$, and for all $A \in \R^{n\times n}, b \in \R^n$, one has the convergence in probability
\[
    \frac{1}{k}\sum_{j=1}^k \norm{\nabla q_{A,b}(\sigma z_j) - \nabla f(x_t + \sigma z_j)}^2 \xrightarrow[k \to \infty]{\Pr_\infty} \int_{\R^n}\norm{\nabla q_{A,b}(z) - \nabla f(z)}^2  \, \Pr_\sigma(\dd z) \,.
\]
\end{assumption}
The next result considers a model, where the objective function is a sum of a quadratic and a function with bounded first derivative, i.e., a \enquote{target} function superimposed with a \enquote{disturbance}. \cref{thm:lanc} describes conditions for which the quadratic approximation objective of \cref{alg:nlqn} can disregard solutions that are not equal to the underlying target model.
\begin{theorem}[Consistency]\leavevmode
\label{thm:lanc}
\noindent Let $\Pr$ be a probability measure with second moments on $\R^n$ and $\Pr_{\sigma} := \Pr(\sigma^{-1}\cdot)$, where $\sigma \in (0, \infty)$.\\
Further, let $f: \R^n \to \R$ be continuously differentiable with
\begin{itemize}
    \item $f \equiv  r + g$, where $r$ is quadratic and $\norm{\nabla g}$ is uniformly bounded by $M \in (0, \infty)$,
    \item $\nabla f$ square-$\Pr$-integrable,
\end{itemize}
and $q^*: \R^n \to \R$ be quadratic, such that,
\[
    \int_{\R^n} \norm{\nabla (q^*-r)(z) - \nabla (q^*-r)(0) }^2 \, \Pr(\mathrm{d}z) > 0 \,.
\tag*{(CON)}
\]
Then there exists $\sigma^* \in (0, \infty)$, such that, for all $\sigma \geq \sigma^*$ the function $q^*$ is suboptimal, i.e.,
\[
    q^{*} \notin \argmin_{q \text{ quadratic}} \int_{\R^n}\norm{\nabla q(z) - \nabla f(z)}^2  \, \Pr_\sigma(\dd z) \,.
\]
\end{theorem}
\begin{remark}\leavevmode
\begin{enumerate}[label=\roman*.]
    \item
    The condition (CON) of \cref{thm:lanc}
    encodes both a failure of $q^*$ to approximate the correct second derivative of $r$ as well as the sampling distribution $\Pr$ to measure this defect.
    \item
    Among the candidates of the minimization problem of \cref{thm:lanc} that have the same second-derivative as $r$, only the optimal ones also have the same first derivatives if $\int_{\R^n} \nabla g \, \dd\Pr = 0$. For this to hold, (CON) nor the boundedness of $\nabla \norm{g}$ are required.
\end{enumerate}
\end{remark}
\begin{proof}
First, it can be seen that $q \overset{!}{=} r$ has an objective value uniformly bounded in $\sigma$, i.e.,
\begin{align*}
    \int_{\R^n} \norm{\nabla q(z) - \nabla f(z)}^2 \,  \Pr_\sigma(\mathrm{d}z)
    &=
    \int_{\R^n} \norm{\nabla g(z)}^2 \,  \Pr_\sigma(\mathrm{d}z)
\tag*{\small{(definition of $f$, and $q = r$)}}
    \\
    &\leq
    M^2
\tag*{\small{(by $\norm{\nabla g} \leq M$, and $\Pr_\sigma$ normed)}}
\end{align*}
It is shown that there exists $\sigma^* \in (0, \infty)$, such that, $q^*$ exceeds this objective value for all $\sigma \geq \sigma^*$, which implies the result. Observe that $\nabla (q-r)$ is affine, i.e., there exist $L :\equiv \nabla (q-r) - \nabla (q-r)(0) \in \R^{n \times n}$ and $w := \nabla (q-r)(0) \in \R^n$, such that, $\nabla (q-r)(z) = L z + w$ for all $z \in \R^n$.\\
In fact, it can be seen that the objective value of $q^*$ diverges to $\infty$ in $\sigma$. One has
\begingroup
\allowdisplaybreaks
\begin{align*}
    &\int_{\R^n} \norm{\nabla q(z) - \nabla f(z)}^2 \,  \Pr_\sigma(\mathrm{d}z)
    \\
    &\qquad\qquad=
    \int_{\R^n} \norm{\nabla (q-r)(z) - \nabla g(z)}^2 \,  \Pr_\sigma(\mathrm{d}z)
\tag*{\small{(definition of $f$, linearity of $\nabla$)}}
    \\
    &\qquad\qquad\geq
    \int_{\R^n} (\norm{\nabla (q-r)(z)} - \norm{\nabla g(z)})^2 \,  \Pr_\sigma(\mathrm{d}z)
\tag*{\small{(the reverse triangle inequality)}}
    \\
    &\qquad\qquad=
    \int_{\R^n} \norm{\nabla (q-r)(z)}^2 - 2\norm{\nabla (q-r)(z)}\norm{\nabla g(z)} + \norm{\nabla g(z)}^2 \,  \Pr_\sigma(\mathrm{d}z)
\tag*{\small{(distributivity)}}
    \\
    &\qquad\qquad\geq
    \int_{\R^n} \norm{\nabla (q-r)(z)}^2 - 2\norm{\nabla (q-r)(z)}\norm{\nabla g(z)} \,  \Pr_\sigma
\tag*{\small{(by $\norm{\cdot} \geq 0$)}}
    \\
    &\qquad\qquad=
    \int_{\R^n} \norm{\nabla (q-r)(z)}^2  \,  \Pr_\sigma(\mathrm{d}z) - 2\int_{\R^n}\norm{\nabla (q-r)(z)}\norm{\nabla g(z)} \,  \Pr_\sigma(\mathrm{d}z)
\tag*{\small{(linearity of $\int$)}}
    \\
    &\qquad\qquad\geq
    \int_{\R^n} \norm{\nabla (q-r)(z)}^2  \,  \Pr_\sigma(\mathrm{d}z) - 2M\int_{\R^n}\norm{\nabla (q-r)(z)} \,  \Pr_\sigma(\mathrm{d}z)
\tag*{\small{(by $\norm{g} \leq M$)}}
    \\
    &\qquad\qquad=
    \int_{\R^n} \norm{L z + w}^2  \,  \Pr_\sigma(\mathrm{d}z) - 2M\int_{\R^n}\norm{Lz + w} \,  \Pr_\sigma(\mathrm{d}z)
\tag*{\small{(by $\nabla (q-r)$ affine)}}
    \\
    &\qquad\qquad=
    \int_{\R^n} \norm{L \sigma z + w}^2  \,  \Pr_\sigma(\mathrm{d}z) - 2M\int_{\R^n}\norm{L \sigma z + w} \,  \Pr(\mathrm{d}z)
\tag*{\small{(by $\Pr_\sigma := \Pr(\sigma^{-1}\cdot)$)}}
    \\
    &\qquad\qquad\geq
    \int_{\R^n} (\norm{L \sigma z} - \norm{w})^2  \,  \Pr(\mathrm{d}z) - 2M\int_{\R^n}\norm{L \sigma z} + \norm{w} \,  \Pr(\mathrm{d}z)
\tag*{\small{(the (reverse) triangle inequality)}}
    \\
    &\qquad\qquad\geq
    \sigma^2 \int_{\R^n} \norm{L z}^2 \,  \Pr(\mathrm{d}z) - 2\sigma (\norm{w}+M)\int_{\R^n} \norm{Lz} \,  \Pr(\mathrm{d}z) + \norm{w}^2 - 2M \norm{w} 
\tag*{\small{(linearity and rearranging)}}
    \\
    &\qquad\qquad\xrightarrow{\sigma \to \infty}
    \infty \,.
\tag*{\small{(by $L z = \nabla (q-r)(z) - \nabla (q-r)(0)$, and by the condition (CON))}}
\end{align*}\qedhere
\endgroup
\end{proof}
The analysis will be continued with an objective function model that is similar to that of \cref{thm:lanc}: A quadratic superimposed by a function which is understood as a \enquote{disturbance}.
The error incurred when selecting the correct quadratic, i.e., the residual, and its dependency of $\sigma$ can generate first insights into which finite values of the scaling $\sigma$ of the non-local kernel are needed. Further, the residual may be informative about which error a good (or even perfect) non-local quadratic approximation of the objective function should have and could be optimized in the design of advanced sampling, scaling, or linesearch methods.\\
A Rastrigin-type model for objective functions has yielded insightful results in the analysis of non-local optimization algorithms in \cite[Theorem~5.2]{mueller} and \cite{schonenberger,omeradzic}. Therefore, in \cref{thm:resb}, a relatively concrete bound for the residual of an optimal approximant in the Rastrigin-type setting is determined, when the samples are independently drawn from a multivariate normal distribution. The result highlights the role of amplitude modulations for this model class.
\begin{theorem}[Residual Bound for a Rastrigin Model]
\label{thm:resb}
    Let $f\equiv r + g:\R^n \to \R$, where
\begin{itemize}
    \item $r$ is a convex quadratic,
    \item $g(x) =  \sum_{j=1}^m a_j \cos(\langle s_j, x \rangle + \psi_j)$ for all $x \in \R^n$, with
    $s_1, \dots, s_m\in \R^n$, $a \in \R^m$, $\psi_1, \dots, \psi_m \in \R$,
    \item such that for some symmetric, positive definite $\Sigma \in \R^{n \times n}$, one has
    \[
        0 < \varepsilon := \min_{\substack{j,\ell \in \N_{\leq m}\\ j\neq \ell}} \{ \norm{s_j + s_\ell}_{\Sigma} , \norm{s_j - s_\ell}_{\Sigma}\} \,.
    \]%
\end{itemize}%
Setting $q \overset{!}{\equiv} r$ and $\Pr_\sigma \overset{!}{=} \mathcal{N}(0, \sigma^2 \Sigma)$, and $S:=(s_1, \dots, s_m) \in \R^{n \times m}$, where $\sigma > 0$, one has
\[
   \int_{\R^n}  \norm{\nabla q(x) - \nabla f(x)}^2 \, \Pr_\sigma(\dd x) \leq  2\norm{a}^2 \norm{S}_{\mathrm{F}}^2 \big(1 + 
     (m-1)\exp( - \sigma^2 \varepsilon^2 /2)\big) \,.
\]
\end{theorem}
\begin{proof}
One has
\begingroup
\allowdisplaybreaks
\begin{align*}
    &\int_{\R^n}  \norm{\nabla q(x) - \nabla f(x)}^2 \, \Pr_\sigma(\dd x)
    \\
    &\qquad\qquad=
    \int_{\R^n}  \norm{\nabla g(x)}^2 \, \Pr_\sigma(\dd x)
\tag*{\small{(by $q \equiv r$ and $f \equiv r + g$)}}
    \\
    &\qquad\qquad=
    \int_{\R^n}  \big\|-\sum_{j=1}^m a_j \sin(\langle s_j, x \rangle + \psi_j) s_j\big\|^2 \, \Pr_\sigma(\dd x)
\tag*{\small{(derivative of $g$)}}
    \\
    &\qquad\qquad=
    \sum_{j=1}^m\sum_{\ell=1}^m a_j a_\ell s_j^T s_\ell \int_{\R^n} \sin(\langle s_j, x \rangle + \psi_j) \sin(\langle s_\ell, x \rangle + \psi_\ell)  \, \Pr_\sigma(\dd x)
    \\
    &\qquad\qquad=
    \frac{1}{4} \sum_{j=1}^m\sum_{\ell=1}^m a_j a_\ell s_j^T s_\ell \bigg(- \exp( -i\psi_j-i\psi_\ell)\int_{\R^n}\exp\big(i\langle -s_j-s_\ell,x \rangle\big) \, \Pr_\sigma(\dd x)
    \\
    &\qquad\qquad\qquad\qquad
    +  \exp(i\psi_j-i\psi_\ell)
    \int_{\R^n}  \exp\big(i\langle s_j-s_\ell,x \rangle \big) \, \Pr_\sigma(\dd x)
    \\
    &\qquad\qquad\qquad\qquad
    +
    \exp(-i\psi_j+i\psi_\ell)\int_{\R^n} \exp\big(i\langle -s_j+s_\ell,x \rangle \big) \, \Pr_\sigma(\dd x)
    \\
    &\qquad\qquad\qquad\qquad
    - \exp(i\psi_j+i\psi_\ell)
    \int_{\R^n}\exp\big(i\langle s_j+s_\ell,x \rangle \big) \, \Pr_\sigma(\dd x)
    \bigg)
\tag*{\small{(Euler formula, linearity)}}
    \\
    &\qquad\qquad=
    \frac{1}{4} \sum_{j=1}^m\sum_{\ell=1}^m a_j a_\ell s_j^T s_\ell \bigg(- \exp( -i\psi_j-i\psi_\ell)\varphi_\sigma(-s_j-s_\ell)
    +  \exp(i\psi_j-i\psi_\ell)\varphi_\sigma(s_j-s_\ell)
    \\
    &\qquad\qquad\qquad\qquad
    +
    \exp(-i\psi_j+i\psi_\ell)\varphi_\sigma(-s_j+s_\ell)
    - \exp(i\psi_j+i\psi_\ell)
    \varphi_\sigma(s_j+s_\ell)
    \bigg)
\tag*{\small{(characteristic function $\varphi_\sigma$ of $\Pr_\sigma$)}}
    \\
    &\qquad\qquad=
    \sum_{j=1}^m\sum_{\ell=1}^m a_j a_\ell s_j^T s_\ell \bigg(
    \cos(\psi_j-\psi_\ell)\exp\big( - \tfrac{1}{2} (s_j - s_\ell)^T \sigma^2 \Sigma (s_j - s_\ell) \big)
    \\
    &\qquad\qquad\qquad\qquad
    -\cos(\psi_j + \psi_\ell)
    \exp\big( - \tfrac{1}{2} (s_j + s_\ell)^T \sigma^2 \Sigma (s_j + s_\ell) \big) \bigg)
\tag*{\small{(by $\Pr_\sigma = \N(0, \sigma^2 \Sigma)$)}}
    \\
    &\qquad\qquad\leq
    2\sum_{j=1}^m a_j^2 \norm{s_j}^2 + 
    2\sum_{\substack{j,\ell=1\\j<\ell}}^m \abs{a_j a_\ell} \abs{s_j^T s_\ell} \bigg(
    \exp\big( - \tfrac{1}{2} (s_j - s_\ell)^T \sigma^2 \Sigma (s_j - s_\ell) \big)
    \\
    &\qquad\qquad\qquad\qquad
    +
    \exp\big( - \tfrac{1}{2} (s_j + s_\ell)^T \sigma^2 \Sigma (s_j + s_\ell) \big) \bigg)
    \tag*{\small{(by $\Delta$-inequality and $\abs{\cos} \leq 1$)}}
    \\
    &\qquad\qquad\leq
    2\sum_{j=1}^m a_j^2 \norm{s_j}^2 + 
    2 \exp( - \sigma^2 \varepsilon^2 /2)\sum_{\substack{j,\ell=1\\j<\ell}}^m 2\abs{a_j a_\ell} \norm{s_j} \norm{s_\ell}
    \tag*{\small{(definition of $\varepsilon$, monotonicity of $\exp$, C-B-S-ineq.)}}
    \\
    &\qquad\qquad\leq
    2\sum_{j=1}^m a_j^2 \norm{s_j}^2 + 
    2 \exp( - \sigma^2 \varepsilon^2 /2)\sum_{\substack{j,\ell=1\\j<\ell}}^m a_j^2 \norm{s_j}^2+ a_\ell^2 \norm{s_\ell}^2
    \tag*{\small{(binomial expansion)}}
    \\
    &\qquad\qquad\leq
    2\sum_{j=1}^m a_j^2 \norm{s_j}^2 \big(1 + 
     (m-1)\exp( - \sigma^2 \varepsilon^2 /2)\big)
    \tag*{\small{(counting summands)}}
    \\
    &\qquad\qquad\leq
    2\norm{a}^2 \norm{S}_{\mathrm{F}}^2 \big(1 + 
     (m-1)\exp( - \sigma^2 \varepsilon^2 /2)\big) \,.
    \tag*{\small{(C-B-S-inequality)}}
\end{align*}\qedhere
\endgroup
\end{proof}
As this work aims to study the basic utility of gradients for global optimization, the formal analysis concludes with the presented conceptual results of \cref{thm:lanc,thm:resb}. The convergence properties of \cref{alg:nlqn} depend on the sampling procedure, sample size $k$, as well as scaling and linesearch techniques, i.e., algorithm parameters and mechanisms that are not the focus of this work. 

\section{Experiments}
\subsection{Comparing the search directions $\Delta x_0$, $-b_0$, and $-\overline{g}$}\label{subsec:exp1}
\begin{figure}[h!]%
\vspace*{10mm}
\begin{center}%
\includegraphics[width=.95\textwidth]{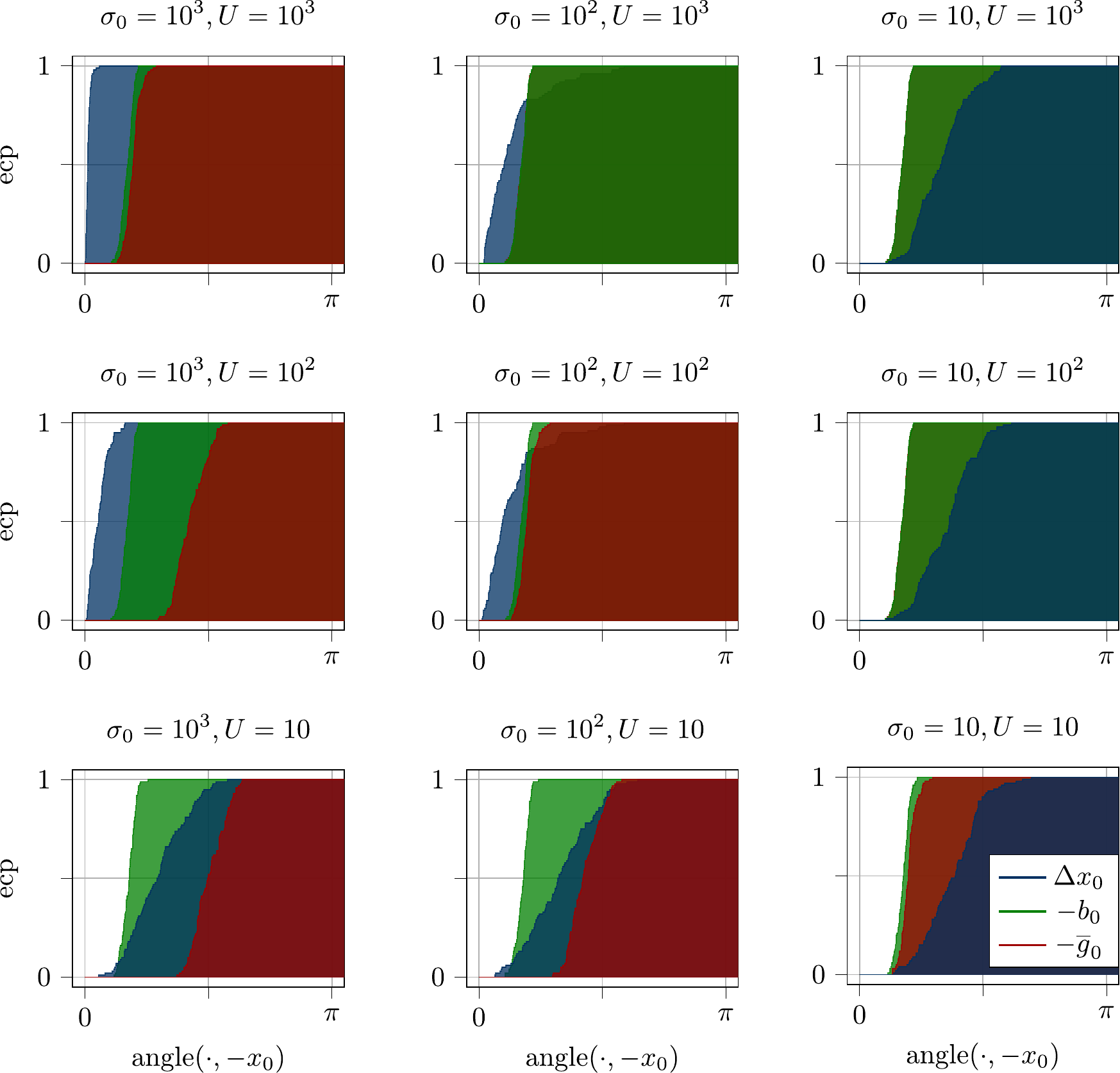}%
\end{center}%
\vspace{3mm}
\caption{Empirical cumulative probability (ecp) over $\mathrm{angle}(\cdot, -x_0) :\equiv \arccos\big(\langle \cdot, -x_0 \rangle / ( \norm{\cdot} \norm{x_0})\big)$ for different estimators of descent directions in the setting of \cref{subsec:exp1} (an ill-conditioned Rastrigin-type function in dimension $n=20$), where $100$ samples of the initial element $(x_0, z_1, \dots, z_k)$ are drawn and $k=30$ gradients are used to estimate the search direction. The vector $-x_0$ is the direction that points towards the global minimum.\\
If the ecp of $\mathrm{angle}(-\overline{g}_0, -x_0)$ is not visible, it is almost equal to that of $\mathrm{angle}(-b_0, -x_0)$.\\
\textbf{Not shown:} Gradient estimators ($-b_0$ and $-\overline{g}$) degrade for $U\leq 1$, whereas $\Delta x_0$ degrades for $U \leq 10$. Further, gradient estimation performs well even for small $\sigma_0$ (i.e., $\sigma_0 = 10^{-2}$), while $\Delta x_0$ degrades for $\sigma_0 \leq 10$. All estimators perform well for the local setting of $\sigma_0, U = 10^2$.}%
\label{fig:descentdirections}%
\end{figure}

The quality of search directions is determined by their deviation from the global optimum. The first experiment measures this deviation with the Euclidean angle between the search direction and the vector that would yield the global optimum (known in the experiment). In a practical setting, however, it should be noted that search directions that deviate consistently by little and those with low deviation in particularly adverse settings are useful---especially when used jointly.\\
In \cref{alg:nlqn}, $\Delta x_0$ denoted the non-local Newton step, $-b_0$ denoted the first order term in the non-local quadratic model, and $-\overline{g}$ denoted in \cref{thm:nsc} the empirical average gradient. For the objective of this experiment, an ill-conditioned quadratic superimposed with a disturbance modeled by several high-frequency cosine functions is chosen. For the selected objective function, second-order information influences the optimal search direction. As in \cref{thm:resb}, this model is an adverse setting that still offers a sufficient degree of interpretation and generality. Therefore, the \textbf{inputs of \cref{alg:nlqn} are set as}
\begin{itemize}
    \item $n \overset{!}{=} 20$, $f(x) \overset{!}{=} \langle x, \mathrm{diag}(1, 1 + (100/19), \dots, 100)x \rangle - \sum_{j=1}^n a \cos( s x_j)$ for all $x \in \R^n$, where $a = 10$ and $s = 20\pi$,
    \item $x_0 \overset{!}{\sim} \mathrm{Unif}([-U, U]^n)$, and where $\sigma_0, U \overset{!}{\in} \{10^{-2}, 10^{-1}, 1, 10, 10^2, 10^3\}$,
    \item $(z_1, \dots, z_k) \sim \Pr_k \overset{!}{=} \NN(0,I_n)^{\otimes k}$, i.e., the product measure of the standard normal distribution of order $k \overset{!}{=} 3n/2$ with $(z_1, \dots, z_k)$ being independent of $x_0$.
    \item The quantities $\Delta x_0, -b_0$ and $-\overline{g}$ do not depend on the other inputs of \cref{alg:nlqn}.
\end{itemize}
The results of this experiment are shown in \cref{fig:descentdirections}. An \textbf{interpretation of the results of the first experiment} is the non-local least-squares gradient estimator $-b_0$ is considerably more robust than the gradient average $-\overline{g}$, as the empirical cumulative probability of $\mathrm{angle}(-b_0, -x_0)$ is larger or almost equal to that of $\mathrm{angle}(-\overline{g}_0, -x_0)$ for all scales of $(\sigma_0, U)$. On large scales of $(\sigma_0, U)$, where the signal-to-noise ratio of gradients of $f$ is large, the non-local Newton direction $\Delta x_0$ is instructive. Here an interpretation is that large $\sigma_0$ improves the estimator, whereas large $U$ makes the estimation problem easier. It is to be noted that estimation is done on relatively few samples $k = 3n/2$. Still, the search directions of \cref{alg:nlqn}, i.e., $\Delta x_0$ and $-b_0$ are clearly significantly more useful than a random search direction in all presented settings. While the experiment studies a particular setting, it is expected that the results generalize whenever the signal-to-noise ratio and the global structure are similar.

\vfill%

\subsection{Benchmarking \cref{alg:nlqn} on selected functions with many local minima}\label{subsec:exp2}
In this experiment, \cref{alg:nlqn}, i.e., the simplistic non-local quasi-Newton method proposed in this work; the \emph{Covariance matrix adaptation evolution strategy} (\texttt{CMA-ES}) \cite{hansen}, an optimizer iteratively modeling up-to-second-moment search distributions from non-local function evaluations; and the \emph{Broyden-Fletcher-Goldfarb-Shanno method} \cite[Chapter~8.1]{nocedal}, a local quasi-Newton method, that is randomly and independently reinitialized uniformly on $[-\sigma_0 + x_0, \sigma_0 + x_0]^n$ on convergence (the algorithm shall be called \texttt{rBFGS}) are compared. The goal is to provide evidence for the utility of non-local gradient information for global optimization in a specific, yet non-trivial and practically relevant, setting. The selected methods \texttt{CMA-ES} and \texttt{rBFGS} are distinct state-of-the-art methods that both model second-order information without explicit access to it (such as \cref{alg:nlqn}).\\
The functions that are used as a benchmark are selected to fit the design domain of \cref{alg:nlqn}, i.e., roughly modeled by a \enquote{disturbed} convex function, and can be considered challenging problems for most optimizers.
In particular, consider the functions defined for all $x \in \R^n$ by
\begin{align*}
    f_{\mathrm{levy}}(x) &:= \sin^2\big(\pi w(x_1)\big) + \big(w(x_n)-1\big)^2 \big(1 + \sin^2\big(2 \pi w(x_n)\big) \big)
    \\
    &\qquad\qquad\qquad\qquad
    + \sum_{i=1}^{n-1} \big(w(x_i) - 1\big)^2\big( 1 + 10 \sin^2\big(\pi w(x_i) + 1\big) \big)
    \,,
    \intertext{\qquad\qquad\qquad\qquad where $w(x_i) := 1+ \frac{x_i - 1}{4}$ for all $i \in \N_{\leq n}$.}
    f_{\mathrm{salomon}}(x) &:= 1 - \cos(12 \pi \norm{x}) + \frac{3}{5} \norm{x} \,, \,\, \text{and}
    \\
    f_{\mathrm{rcigar}}(x) &:= an + \langle x, \mathrm{diag}(1, 1 + 100/(n-1), \dots, 100)x \rangle - \sum_{i=1}^n a \cos( s x_i)
    \,,\,\, \text{where $a=10$ and $s=20\pi$} \,.
\end{align*}
The \textbf{inputs of \cref{alg:nlqn} are set as}
\begin{itemize}
    \item $n \overset{!}{=} 50$, $f \overset{!}{\in} \{f_{\mathrm{levy}}, f_{\mathrm{salomon}}, f_{\mathrm{rcigar}}\}$, $x_0 \overset{!}{\sim} \mathrm{Unif}([-10, 10]^n)$, $\sigma_0 \overset{!}{=} 10$, and
    \item $(z_1, \dots, z_k) \sim \Pr_k \overset{!}{=} \NN(0,I_n)^{\otimes k}$, i.e., the product measure of the standard normal distribution of order $k = 3n$ with $(z_1, \dots, z_k)$ being independent of $x_0$.
    \item Define 
    \begin{align*}
        \mathrm{scaling}(\sigma_0, \sigma_t, x_{t+1}-x_t) &:=
        \begin{cases}
            \mathrm{scaling}(\sigma_0, \sigma_0, x_{t+1}-x_t)
            & \text{if } \sigma_t < 10^{-4}
            \\
            \sigma_t/2
            & \text{else if } \norm{x_{t+1}-x_t} < 10^{-4}
            \\
            \norm{x_{t+1} -x_t}/2
            & \text{else if } \norm{x_{t+1}-x_t} > 2\sigma_t
            \\
            \sigma_t & \text{else,}
        \end{cases}
        \\
        \text{and} \quad \mathrm{linesearch}(f, \Delta x_t, -b_t, x_t) &:=
        \argmin_{x \in \mathcal{A}(\Delta x_t, -b_t, x_t)} f(x) \,,
    \intertext{\qquad\qquad\qquad\qquad where $\mathcal{A}(\Delta x_t, -b_t, x_t) := \{x_t + (6/5)^i \Delta x_t, x_t + (6/5)^i (-b_t) \mid i\in \Z_{[-10, 10]}\}$.}
    \end{align*}
\end{itemize}
The initial scaling of \texttt{CMA-ES}, i.e., the size of the non-local kernel, was chosen to be the same as for \cref{alg:nlqn}. Virtually all other parameters of \texttt{CMA-ES} are heavily tuned for the considered function classes and adaptive to the functions at hand.\\
The local search of \texttt{rBFGS} is terminated (and restarted) when the gradient norm is smaller than $10^{-4}$. The other parameters of \texttt{rBFGS} do not affect its performance significantly, as local search is extremely efficient on the considered functions. The reinitialization will likely not be improved much by gridded sampling as the experiments are done in $50$ dimensions, i.e., reinitializations are very sparse.
\begin{figure}[h!]%
\vspace{3mm}
\begin{center}%
\includegraphics[width=.95\textwidth]{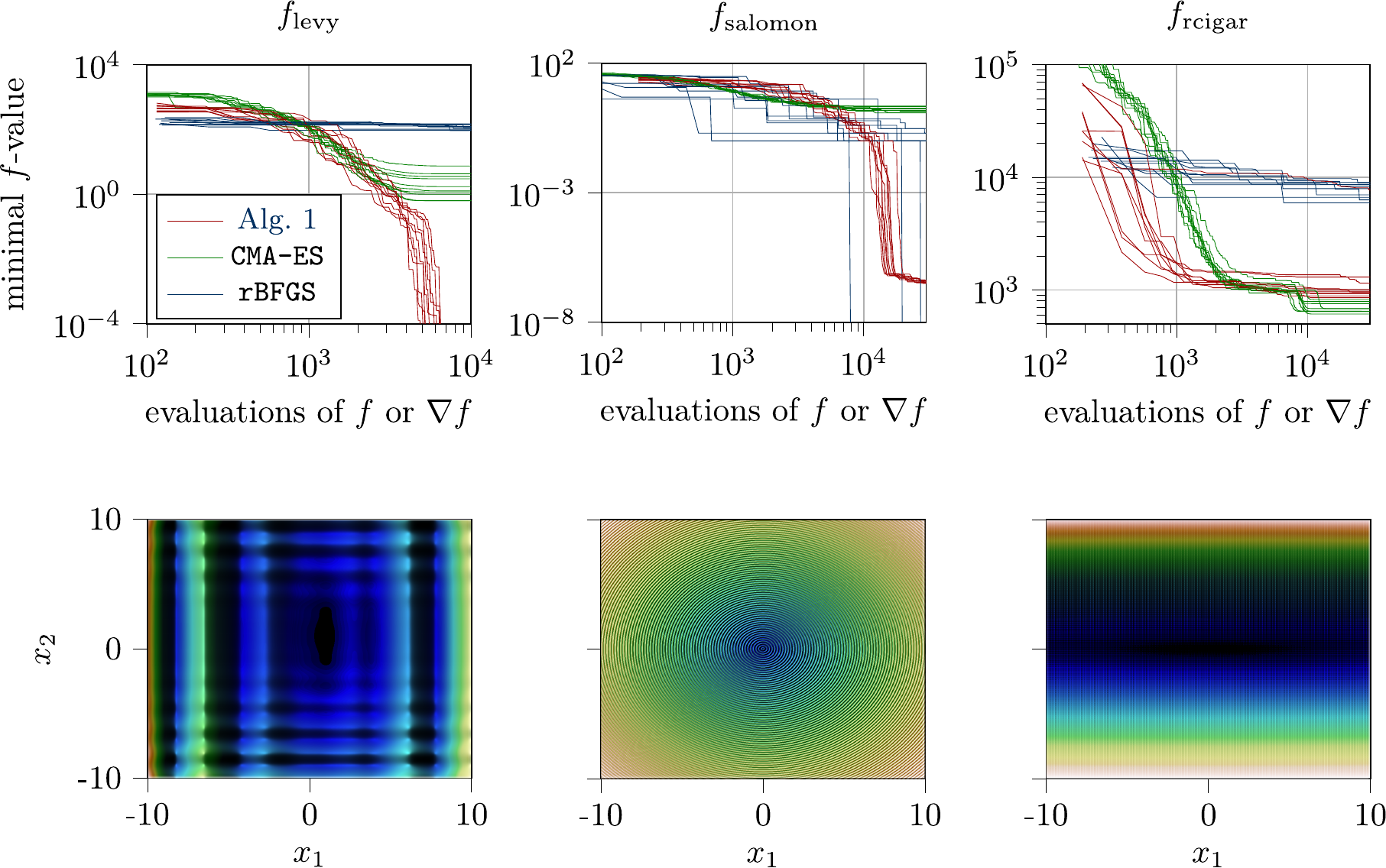}
\end{center}%
\caption{\textbf{Top:} Benchmark of the proposed \cref{alg:nlqn}, zero-order stochastic search method \texttt{CMA-ES} and uniformly randomly reinitialized local quasi-Newton method \texttt{rBFGS}. The benchmark functions $f_{\mathrm{levy}}, f_{\mathrm{salomon}}$ and $f_{\mathrm{rcigar}}$ in dimension $n=50$, defined in \cref{subsec:exp1} are selected for their pathological structure with many local minima. The global minimum of all functions is $0$. Function and gradient evaluations (including the linesearch) are counted equally. All algorithms are initialized independently for each run uniform randomly on $[-10, 10]^n$.\\
\textbf{Bottom:} Visualization of $2$-dimensional analogs of the benchmark functions. Lighter/reddish colors correspond to relatively large function values and darker/bluish colors correspond to relatively small function values.}
\label{fig:benchmark}
\end{figure}
The results of this experiment are shown in \cref{fig:benchmark}.\\
An \textbf{interpretation of the results of the second experiment} is that for functions with benign second-order structure, such as $f_\mathrm{levy}$, \cref{alg:nlqn} and \texttt{CMA-ES} perform similarly well on large-scales, as non-local approximates of up-to-second order objective structure seems to yield good search directions and both methods estimate them. \cref{alg:nlqn} seems to have a slight advantage due to gradients being more informative than function values. In the particular case of $f_\mathrm{levy}$ the likelihood of reinitializing in a relatively good basin of attraction is small, such that, \texttt{rBFGS} performs relatively poorly. It is plausible that due to its robust linesearch, \cref{alg:nlqn} enters a basin of attraction of $f_\mathrm{levy}$ close to the minimum and fast local convergence using second-order information is possible. This is evidence for the hypothesis that non-local methods profit from linesearch. Similar behavior is observed for $f_\mathrm{salomon}$, while here it is likely enough that \texttt{rBFGS} reinitializes in a good basin of attraction given the available budget. Therefore, \texttt{rBFGS} performs relatively well on $f_\mathrm{salomon}$. The reduction in the rate of convergence of \cref{alg:nlqn} of $f_\mathrm{salomon}$ is likely due to its coarse stepping in the scaling function possibly paired with the non-smoothness in the origin of $f_\mathrm{salomon}$. Despite the forgetfulness of \cref{alg:nlqn} (function evaluations are not retained past an iteration), \cref{alg:nlqn} seems to outperform \texttt{CMA-ES} in estimating second-order information in domains of benign signal-to-noise-ratio of $f_\mathrm{rcigar}$. \texttt{CMA-ES} on the contrary seems to perform relatively better in the domain of poor signal-to-noise-ratio close to the optimum of $f_\mathrm{rcigar}$. The many relatively bad local minima of $f_\mathrm{rcigar}$ yield a relatively bad performance of \texttt{rBFGS}.

\subsection{Solving Problem 4 of \emph{SIAM News:~A Hundred-dollar, Hundred-digit Challenge}}\label{subsec:exp3}

\begin{figure}[h!]%
\vspace*{7mm}
\begin{center}%
\includegraphics[width=.95\textwidth]{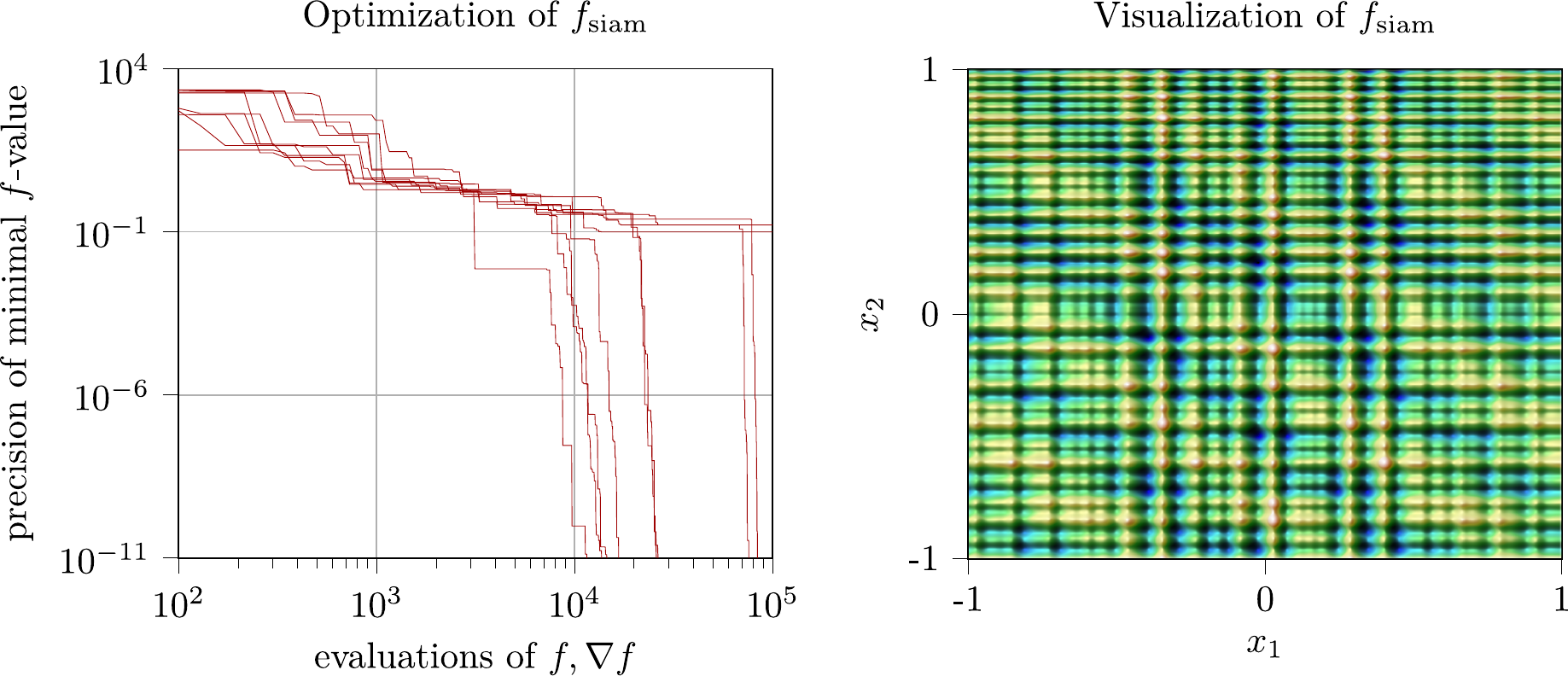}%
\end{center}%
\caption{\textbf{Left:} Benchmark of the proposed \cref{alg:nlqn} on Problem 4 of the \emph{Hundred-dollar, Hundred-digit Challenge}, posed on January 2, 2002 in SIAM News \cite{trefethen}. The algorithm is initialized independently for each run uniform randomly on $[-100, 100]^2$. Function and gradient evaluations (including the linesearch) are counted equally.\\
\textbf{Right:} Visualization of $f_{\mathrm{siam}}$ on $[-1, 1]^2$, which contains the global minimum $-3.306868647475\dots$ at $x^* = (-0.0244\dots, 0.2106\dots)$. Lighter/reddish colors correspond to relatively large function values and darker/bluish colors correspond to relatively small function values.}
\label{fig:siam}
\end{figure}

The \emph{Hundred-dollar, Hundred-digit Challenge}, posed on January 2, 2002 in SIAM News \cite{trefethen}, asks to solve $10$ numerical problems with a precision of $10$ digits each. Problem 4 of the challenge is the minimization of
\begin{align*}
    &x \in \R^2 \longmapsto f_{\mathrm{siam}}(x) := \exp\big(\sin(50x_1)\big) + \sin\big(60 \exp(x_2)\big) + \sin\big(70 \sin(x_1) \big)\\
    &\qquad\qquad\qquad\qquad\qquad\qquad
    + \sin\big( \sin(80x_2) \big) - \sin\big( 10 (x_1 + x_2) \big) + (x_1^2 + x_2^2)/4 \,.
\end{align*}
Although many fundamentally different solutions to Problem 4 exist \cite{bornemann}, it was proposed as a particularly challenging example for many optimization methods. Therefore, it is interesting to verify whether the newly proposed \cref{alg:nlqn} can solve this familiar challenge. When uniform randomly initialized on $[-100,100]^2$, \cref{alg:nlqn} solves the problem in many cases within $30000$ evaluations. The results of the experiment are visualized in \cref{fig:siam}.

For this experiment the \textbf{inputs of \cref{alg:nlqn} were set as}
\begin{itemize}
    \item $n=2, f \overset{!}{=} f_{\mathrm{siam}}$, $x_0 \overset{!}{\sim} \mathrm{Unif}([-100, 100]^n)$, $\sigma_0 \overset{!}{=} 1$, and
    \item $(z_1, \dots, z_k) \sim \Pr_k \overset{!}{=} \NN(0,I_n)^{\otimes k}$, i.e., the product measure of the standard normal distribution of order $k = 3$ with $(z_1, \dots, z_k)$ being independent of $x_0$.
    \item Define 
    \begin{align*}
        \mathrm{scaling}(\sigma_0, \sigma_t, x_{t+1}-x_t) &:=
        \begin{cases}
            \mathrm{scaling}(\sigma_0, \sigma_0, x_{t+1}-x_t)
            & \text{if } \sigma_t < 10^{-4}
            \\
            10\sigma_t/11
            & \text{else if } \norm{x_{t+1}-x_t} < 10^{-4}
            \\
            10\norm{x_{t+1} -x_t}/11
            & \text{else if } \norm{x_{t+1}-x_t} > 2\sigma_t
            \\
            \sigma_t & \text{else.}
        \end{cases}
    \end{align*}
\end{itemize}
All other elements of \cref{alg:nlqn} are defined as in \cref{subsec:exp2}.

\section{Discussion and Future Work}\label{sec:discussion}
This work asked the question of whether practically realistic numbers of gradient evaluations can yield useful information for minimizing functions with many suboptimal local minima.\\
The proposed simplistic \cref{alg:nlqn} generalizes the quasi-Newton method based on a non-local approximation of the objective function. When optimizing a \enquote{disturbed} quadratic function, \cref{alg:nlqn} is a consistent approximator of the underlying quadratic model in the sense of \cref{thm:lanc}. Further, at least one of the two generated search directions outperforms gradient averaging in \cref{subsec:exp1}.\\
Similarly, even with simple sampling, scaling, or linesearch methods, \cref{alg:nlqn} performs comparable to or outperforms state-of-the-art methods on minimization of differentiable functions with many suboptimal local minima selected in \cref{subsec:exp2}. An additional motivation for the extension of the presented method is that \cref{alg:nlqn} solves Problem 4 of the \emph{SIAM News: A Hundred-dollar, Hundred-digit Challenge} \cite{trefethen,bornemann}---likely as (one of) the first non-global non-local iterative methods with search directions based on gradients.\\
The results of this work motivate the design/analysis/integration of sophisticated sampling, scaling, and linesearch methods that adapt to the particular objective function at hand. More precisely, promising future work builds on the described results by
\begin{itemize}
    \item Regularizing and improving the approximation of the quadratic model based on objective evaluations;
    \item Designing improved $\mathrm{linesearch}$ and $\mathrm{scaling}$ methods towards an objective-adaptive optimizer;
    \item Integrating objective and objective gradient evaluations of multiple algorithm iterations to improve the sample efficiency.
\end{itemize}

\newpage

\printbibliography
\end{document}